\documentclass[12pt]{amsart}
\usepackage[cp1251]{inputenc}
\usepackage[english]{babel}

\usepackage{amssymb,amsmath,amscd,amsthm}

\def\{{\protect\lbrace}
\def\}{\protect\rbrace}

\newcommand{\Kdim}{\operatorname{Kdim}}

\begin{document}
\begin{center}
\large{\textbf{Distributive Noetherian Centrally Essential Rings}}
\end{center}

\hfill {\sf Victor Markov}

\hfill Lomonosov Moscow State University

\hfill e-mail: vtmarkov@yandex.ru

\hfill {\sf Askar Tuganbaev}

\hfill National Research University "MPEI"

\hfill Lomonosov Moscow State University

\hfill e-mail: tuganbaev@gmail.com

{\bf Abstract.} It is proved that a ring $A$ is a right or left Noetherian, right distributive centrally essential ring if and only if $A=A_1\times\cdots\times A_n$, where each of the rings $A_i$ is either a commutative Dedekind domain or a  uniserial Artinian centrally essential (not necessarily commutative) ring.

V.T.Markov is supported by the Russian Foundation for Basic Research, project 17-01-00895-A. A.A.Tuganbaev is supported by Russian Scientific Foundation, project 16-11-10013.

{\bf Key words:} centrally essential ring, distributive ring, uniserial ring, Noetherian ring, Artinian ring.

\textbf{MSC2010 datebase 16D25; 16D10}

\section{Introduction}
All rings considered are associative and contain a non-zero identity element. Writing expressions of the form ``$A$ is an Artinian ring'' mean that the both  modules $A_A$ and $_AA$ are Artinian.

A ring $A$ is said to be \emph{centrally essential}\footnote{Centrally essential rings are studied, for example, in \cite{MT18}, \cite{MT19}, \cite{MT20}.} if for any non-zero element $a\in A$, there are two non-zero central elements $c,d\in A$ with $ac=d$, i.e., the module $A_C$ is an essential extension of the module $C_C$, where $C$ is the center of $A$. 

It is clear that all commutative rings are centrally essential. If $Z_2$ is the field of order 2 and $Q_8$ is the quaternion group of order 8, then the group ring $Z_2[Q_8]$ is an example of a non-commutative centrally essential ring \cite{MT18}; in addition, $Z_2[Q_8]$ is a finite local ring. 

A module is said to be \textit{distributive} if the lattice of all its submodules is distributive. A module is said to be \textit{uniserial} if any two its submodules are comparable with respect to inclusion. It is clear that every uniserial module is distributive. The ring of integers $\mathbb{Z}$ is an example of a commutative distributive non-uniserial ring.

This paper is a continuation of \cite{MT20} which is a continuation of \cite{MT19}. In Theorem 1.1, we recall some results of these papers.

\textbf{Theorem 1.1; \cite{MT20}.}~ 
\textbf{1)} Any finite, right uniserial, centrally essential ring is commutative and there are non-commutative uniserial Artinian centrally essential rings; see \cite{MT19}.

\textbf{2)} A centrally essential ring $R$ is a right uniserial, right Noetherian ring if and only if $R$ is a commutative discrete valuation domain or a (not necessarily commutative) uniserial Artinian ring.; see \cite{MT20}.

In connection with Theorem 1.1, we prove Theorem 1.2 which is the main result of this paper.

\textbf{Theorem 1.2.}~
\textbf{1)} If $A$ is a right distributive centrally essential finite ring, then  $A=A_1\times \cdots \times A_n$, where $A_i$ is a commutative uniserial finite ring, $i=1,\ldots ,n$. In particular, $A$  is a commutative distributive ring.\\
Conversely, any such a direct product is a commutative distributive finite ring.

\textbf{2)} If $A$ is a right distributive, right or left Noetherian, centrally essential ring, then  $A=A_1\times \cdots \times A_n$, where $A_i$ is either a (not necessarily commutative) uniserial Artinian ring or a commutative Dedekind domain, $i=1,\ldots ,n$. In particular, $A$  is a distributive Noetherian ring.\\ Conversely, any such a direct product is a distributive Noetherian ring which is not necessarily commutative, since there exist non-commutative uniserial Artinian centrally essential rings; see Theorem 1.1(1).

The proof of the Theorem 1.2 is given in the next section and is based on several assertions, some of which are of independent interest.

We give some notation and definitions. 

Let $M$ be a module. The module $M$ is said to be \textit{uniform} if any two its non-zero submodules have the non-zero intersection. A module $M$ is called an \textit{essential extension} of some its submodule $X$ if $X\cap Y\ne 0$ for any non-zero submodule $Y$ in $M$. For any ring $A$, we denote by $J(A)$, $A^*$ and $C(A)$ the Jacobson radical, the group of invertible elements and the center of  $A$, respectively.

The \emph{left annihilator} of an arbitrary subset $S$ of the ring $A$ is $\ell_A(S)=\{a\in A\,|\, aS=0\}.$
The right annihilator $r(S)$ of the set $S$ are similarly defined.
We denote by $[a,b]=ab-ba$ for the \emph{commutator} of the elements $a,b$ of an arbitrary ring.

Let $A$ be a ring and $B$ a proper ideal of $A$. A ring $A$ is called a \textit{domain} if $A$ does not have non-zero zero-divisors. The ideal $B$ is said to be \textit{completely prime} if the factor ring $A/B$ is a domain.

Other ring-theoretical notions and notations can be found in \cite{Herstein, Lambek, Tug98}.

\section{The proof of Theorem 1.2}\label{section2}

\textbf{Lemma 2.1.} Let $A$ be a centrally essential ring such that the set $B$  of all left zero-divisors of the ring $A$ is an ideal of the ring $A$. Then $B$ coincides with the set $B'$ of all right zero-divisors of the ring $A$ and $A/B$~is a commutative domain. In addition, if the ideal $B$ is nilpotent of nilpotence index $n+1$, then ideal $B^n$ is contained in the center of the ring $A$.

\textbf{Proof.} Let $a,a'\in A\setminus B$. 

If $aa'\in B$, then there exists a non-zero element $x\in A$ with $aa'x=0$. Since $a'\notin B$, we have $a'x\ne 0$. Since $a\notin B$,  we have $aa'x\ne 0$. This is a contradiction. Therefore, $A/B$ is a domain.

Since $A$ is a centrally essential ring, there exist two non-zero central elements $x,y$ of $A$ with  $b x=y$. Then $[a,b]x=[a,bx]=0$, i.e., $[a,b]$~is a left zero-divisor. Consequently, $[a,b]\in B$ and ring $A/B$ is commutative.

Let  $a_1,a_2$ be two non-zero elements of the centrally essential ring $A$ and  $a_1a_2=0$. There exist four non-zero central elements $x_1,x_2,y_1,y_2$ of $A$ such that  $a_1 x_1=y_1$, $a_2 x_2=y_2$. Then $y_2a_1 =a_1y_2 =a_1a_2x_2=0$. Therefore, the left zero-divisor $a_1$~is a right zero-divisor. Similarly, the right zero-divisor $a_2$~is a left zero-divisor. Therefore, $B=B'$.

Now let the ideal $B$ be nilpotent of nilpotence index $n+1$. It remains to prove that $ab-ba=0$ for any elements $a\in A$ and $b\in B^n$. 

Let's assume that $ab-ba\ne 0$ for some elements $a\in A$ and $b\in B^n$. Since $A$ is a centrally essential ring, there exist central elements $x,y\in A$ with $bx=y\ne 0$. Then $abx=ay=ya=bxa=bax$, $(ab-ba)x=abx-bax=0$. Therefore, $x$ is a central zero-divisor. Consequently, $x\in B$. Then $y=bx\in B^{n+1}=0$. This is a contradiction.~\hfill$\square$

For convenience, we give the proof of the following well known Lemma 2.2.

\textbf{Lemma 2.2.}\label{div-fg-mod} Let $A$ be a commutative domain and  there exists a non-zero finitely generated divisible torsion-free $A$-module $M$. Then $A$~is a field.

\textbf{Proof.} We assume the contrary. Then $A$ has a maximal ideal $\mathfrak m\neq 0$ and $M$ can be naturally turn into a non-zero finitely generated module over the local ring $A_{\mathfrak m}$ with radical $J={\mathfrak m}A_{\mathfrak m}$. Since the module $M$ is divisible, we have that
$MJ\supseteq M{\mathfrak m}=M$, and $M=0$ by the Nakayama lemma. This is a contradiction.~\hfill$\square$

\textbf{Proposition 2.3.} Let $A$ be a right distributive indecomposable ring with the maximum condition on right annihilators, which does not contain an infinite direct sum of non-zero right ideals, and let $B$ be the prime radical of the ring $A$. Suppose that $A$ be not a right uniform domain. Then:

\textbf{2.3.1}. The ideal $B$ is a non-zero completely prime nilpotent ideal, $B$ is the set of all left zero-divisors of the ring $A$, $B$ is a essential right ideal and $B=aB$ for any element $a\in A\setminus P$.

\textbf{2.3.2}. If the ring $A$ is centrally essential and $B$ is a finitely generated right ideal of the ring $A$, then  $A$ is a uniserial Artinian ring with radical $B$.

\textbf{Proof.} \textbf{2.3.1}. The assertion is proved in \cite[The assertion 9.20(1)]{Tug98}.

\textbf{2.3.2}. By 2.2.1 and Lemma 2.1, $B$ is a non-zero completely prime nilpotent ideal, $B$ coincides with the set of all left or right zero-divisors of the ring $A$, $A/B$~is a commutative domain and $B=aB$ for any element $a\in A\setminus B$. 

Let $n+1$ be the of nilpotence index of the ideal $B$. Then $B^n$ is a left and right module over the commutative domain $A/B$. Let $a\in A\setminus B$. Since $B=aB$ and $B$ is an ideal, $B^n=aB^n$ and $B^n$ is a divisible left $A/B$-module.  By Lemma 2.1, the ideal $B^n$ is contained in the center of the ring $A$. Therefore, the divisible left $A/B$-module $B^n$ is a divisible right $A/B$-module. Since $B$ is a finitely generated right ideal which is an ideal, $B^n$ is a finitely generated right ideal. Therefore, the commutative domain $A/B$ has a non-zero finitely generated divisible right module $B^n$. By Lemma 2.1, any element of the set $A\setminus B$ is not a right zero-divisor. Therefore, $B^n$  is a torsion-free right $A/B$-module. By Lemma 2.2, $A/B$ is a field. Therefore, $A$ is a right distributive local ring with non-zero nilpotent Jacobson radical $B$. Then $A$ is a right uniserial Artinian right ring \cite{Ste74}. Since $A$ is a right uniserial, right Artinian, centrally essential ring, $A$ is a uniserial Artinian ring \cite[Corollary 2.4]{MT19}.~\hfill$\square$

\textbf{Corollary 2.4.} Let $A$ be a right distributive Noetherian right centrally essential indecomposable ring with non-zero prime radical $B$. Then $A$ is a uniserial Artinian ring with radical $B$.

Corollary 2.4 follows from Lemma 2.3.2.
 
\textbf{Remark 2.5.} A semiprime ring $A$ is centrally essential if and only if $A$ is commutative \cite[Proposition 3.3]{MT18}. Consequently, a prime  ring is a centrally essential if and only if $A$ is a commutative domain. 

\textbf{Remark 2.6.} It is well known that commutative domain $A$ is a distributive Noetherian ring if and only if $A$ is a Dedekind domain; e.g., see \cite[Theorem 6]{Cam75} or \cite[Theorem 9.16]{Tug98}. It follows from this property and Remark 2.5 that the right distributive, right or left Noetherian, centrally essential prime rings coincide with commutative Dedekind domains. Such rings are distributive Noetherian rings.

\textbf{Proposition 2.7.} A ring $A$ is a right distributive,  right Noetherian centrally essential ring if and only if $A=A_1\times \cdots \times A_n$, where $A_i$ is a uniserial Artinian centrally essential ring or commutative Dedekind domain, $i=1,\ldots ,n$. Conversely, it is clear that in this case $A$ is a distributive Noetherian ring.

Proposition 2.7 follows from Corollary 2.4 and Remarks 2.5 and 2.6.

\textbf{Proposition 2.8.} A ring $A$ is a right distributive, left Noetherian centrally essential ring if and only if $A=A_1\times \cdots \times A_n$, where $A_i$ is a uniserial Artinian centrally essential ring or a commutative Dedekind domain, $i=1,\ldots ,n$. Conversely, it is clear that $A$ is a distributive Noetherian ring in this case.

\textbf{Proof.} In \cite[Theorem 9.18(1)]{Tug98}, it is proved that a ring $A$ is a right distributive, left Noetherian ring if and only if $A=A_1\times \cdots \times A_n$, where $A_i$ is an Artinian, right uniserial ring or right distributive, left Noetherian domain, $i=1,\ldots ,n$. In addition, any right uniserial centrally essential Artinian ring is a left uniserial ring \cite[Corollary 2.4]{MT19}. Now we use Remark 2.6.~\hfill$\square$ 

\textbf{2.9. The completion of the proof of Theorem 1.2.}\\
\textbf{1.2(1).} The right distributive finite rings coincide with finite direct products of right uniserial finite rings \cite{Ste74}. Therefore, 1.2(1) follows from Theorem 1.1(1).\\
\textbf{1.2(2).} The assertion follows from Propositions 2.7 and 2.8.

\section{Remarks}
\mbox{~}

\textbf{3.1. Remark.} In \cite{Jen64}, it is proved that a commutative ring $A$ is a distributive semiprime ring if and only if every submodule of any flat\footnote{A right $A$-module $X$ is said to be \textsf{flat} if for every left $A$-module $Y$ and any submodule $Y'$ in $Y$, a natural homomorphism $X\otimes_AY'\to X\otimes_AY$ of additive groups is a monomorphism.} $A$-module is a flat module. Therefore, it follows from Remark 2.5 that the right or left distributive, centrally essential, semiprime rings coincide with commutative rings over which every submodule of any flat module is a flat module.

\textbf{3.2. Remark.} In \cite{Cam75}, it is proved that the commutative distributive semiprime rings, not containing an infinite direct sum of non-zero ideals, coincide with finite direct products of commutative semihereditary\footnote{A module is said to be \textsf{semihereditary} if all its finitely generated submodules are projective modules.} domains. Therefore, it follows from Remark 2.5 that the right or left distributive, centrally essential, semiprime rings, not containing an infinite direct sum of non-zero ideals, coincide with finite direct products of commutative semihereditary domains.

\textbf{3.3. Modules with Krull dimension.} We recall the transfinite definition of the Krull dimension $\Kdim M$ of the module $M$, see  \cite{GorR73}. \\ 
By definition, we assume that zero modules are of Krull dimension $-1$ and every non-zero Artinian module is of Krull dimension $0$.\\ 
Let's assume that $\alpha$ is an ordinal number $>0$ such that the modules with Krull dimension $\beta$ have been defined for all ordinal numbers $\beta <\alpha$ and let $M$ be a module with $\Kdim M\ne \beta$.\\ 
We say that the \textsf{Krull dimension}\label{kdim} $\Kdim M$ of the module $M$ is $\alpha$ if for any infinite strictly descending chain $M_1>M_2>\ldots$ of submodules of $M$, there exists a positive integer 
$n$ with $\Kdim (M_n/M_{n+1})<\alpha$. 

If $A$ is a ring and the Krull dimension $\Kdim A_A$ of the module $A_A$ exists then $\Kdim A_A$ is called the \textsf{right Krull dimension} of the ring $A$.

\textbf{3.3.1.} Every Noetherian module has Krull dimension; see \cite{GorR73}.

\textbf{3.3.2.} Not every module has Krull dimension; for example, any direct sum of an infinite set of non-zero modules does not have Krull dimension.

\textbf{3.4. An example of a commutative, uniserial, non-Artinian, non-prime ring of Krull dimension 1.}\\
It follows from Theorem 1.2(2) that every right distributive, right Noetherian, centrally essential, indecomposable ring $A$ is either a commutative domain or an Artinian ring. In this assertion, we cannot replace the condition "$A$ is a right Noetherian ring" by the condition "$A$ is a ring with right Krull dimension". Indeed, let $D$ be a commutative uniserial principal ideal domain which is not a field and let $Q$ be the field of fractions of the domain $D$. For example, we can take the formal power series ring $D=F[[x]]$ over a field $F$ and the formal Laurent series field $Q=F((x))$. Let $A$ be the trivial extension of the $D$-$D$-bimodule $Q$ by $D$. It is directly verified that $A$ is a commutative, uniserial, non-Artinian, non-prime ring of Krull dimension 1.

\textbf{3.5. An example of a commutative, distributive, non-uniserial, non-Artinian, non-prime ring of Krull dimension 1.}\\
It follows from Theorem 1.2(2) that every right distributive, right Noetherian, centrally essential, indecomposable ring $A$ is either a commutative domain or a uniserial ring. In this assertion, we cannot replace the condition "$A$ is a right Noetherian ring" by the condition "$A$ is a ring with right Krull dimension". Indeed, let $\mathbb{Z}$ be the ring of integers, $p$ be a prime integer, $M$ be the quasi-cyclic  $p$-group, and let $A$ be the trivial extension of the $\mathbb{Z}$-$\mathbb{Z}$-bimodule $M$ by $\mathbb{Z}$. It is directly verified that $A$ is a commutative distributive, non-uniserial, non-Artinian, non-prime, indecomposable ring of Krull dimension 1. In addition, $A$ is a Bezout ring in which any ideal is comparable with respect to inclusion with the ideal $\mathbb{Q}$.

\end{document}